\documentclass[12pt]{article}%
\usepackage{amsmath}
\usepackage{amsfonts}
\usepackage{amssymb}
\usepackage{graphicx}%
\providecommand{\U}[1]{\protect\rule{.1in}{.1in}}
\newtheorem{theorem}{Theorem}

\newtheorem{corollary}{Corollary}

\newtheorem{example}{Example}

\newtheorem{lemma}{Lemma}

\newtheorem{proposition}{Proposition}

\begin{document}

\title{Poisson limit theorem for the number of excursions above high and medium levels by Gaussian
stationary sequences}
%\author{V. I. Piterbarg}
\author{Vladimir I. Piterbarg\thanks{Lomonosov Moscow state university, Moscow,
Russia; Federal State Institution "Scientific-Research Institute for System Analysis of the Russian Academy of Sciences"; International Laboratory of Stochastic Analysis and its
Applications, National Research University Higher School of Economics, Russian
Federation.\ \emph{{piter@mech.math.msu.su}}}}
\maketitle

\begin{abstract}
Asymptotic behavior of the point process of high and medium values of a
Gaussian stationary process with discrete time is considered. An approximation
by a Poisson cluster point process is given for the point process.

\end{abstract}

Consider Gaussian stationary sequence $X(k),$ $k\in\mathbb{Z}$ with $EX(k)=0,$
$\operatorname*{Var}X(k)\equiv1,$ $\operatorname*{Cov}%
(X(0,X(k)=EX(k)X(0)=r(k).$ Let $\mathcal{B}$ be the algebra of bounded Borel
subsets of $\mathbb{R}$. Consider on $\mathcal{B}$ a family of Bernoulli
processes
\begin{equation}
\mathfrak{B}_{u,n}(B):=\sum_{k\in nB}\mathbf{I}\left\{  X(k)>u\right\}
,\ \ B\in\mathcal{B},\ u>0,\ n\in\mathbb{N}. \label{BPP}%
\end{equation}
We study the limit behavior of $\mathfrak{B}_{u,n}(\cdot)\ $as $u,n\rightarrow
\infty.$ Theorem \ref{Mittal Ylvisaker} below says that if $r(k)$ tends to
zero sufficiently fast, $\mathfrak{B}_{u,n}(\cdot)$ tends as $u,n\rightarrow
\infty$ weakly to a Poisson point process $\mathfrak{P}_{\lambda}(B),$
$B\in\mathcal{B}$, with intensity $\lambda>0$ provided \emph{natural
normalization} is fulfilled, that is, $n,u\rightarrow\infty$ such that
\begin{equation}
\lim_{u,n\rightarrow\infty}np(u)=\lambda,\ p(u)=P(X(1)>u). \label{norm0}%
\end{equation}
In case of independent $X(k)$s it means that the distribution (binomial) of
$\mathfrak{B}_{u,n}([0,1])$ tends to Poisson one with parameter $\lambda$
(Poisson Limit Theorem). Since for Gaussian standard distribution function
$\Phi$,
\begin{equation}
\Psi(u)\geq p(u)=1-\Phi(u)\geq(1-u^{-2})\Psi(u),\ \ u>0,\ \text{with }%
\Psi(u)=\frac{1}{\sqrt{2\pi}u}e^{-u^{2}/2}, \label{p(u)}%
\end{equation}
analytically more convenient to take an equivalent normalization, with
$\Psi(u)$ instead of $p(u).$ From relations (\ref{norm0},\ref{p(u)}) one has
the asymptotic solution,
\begin{equation}
u=u_{n}=\sqrt{2\log n}-\frac{\frac{1}{2}\log\log n+\log(\lambda\sqrt{\pi/2}%
)}{\sqrt{2\log n}}+O(1/\log n),\ n\rightarrow\infty. \label{norm}%
\end{equation}

From results of \cite{mittal ylv} it follows

\begin{theorem}
\label{Mittal Ylvisaker} Let
\begin{equation}
r(k)\log k\rightarrow0\ \text{as }k\rightarrow\infty, \label{logCond}%
\end{equation}
and (\ref{norm0}) be fulfilled. Then
\begin{equation}
\mathfrak{B}_{u,n}(B)\Rightarrow\mathfrak{P}_{\lambda}(B),\ B\in\mathcal{B},
\label{poisson}%
\end{equation}
weakly as $u\rightarrow\infty.$
\end{theorem}

Mittal and Ylvisaker in \cite{mittal ylv} proved a limit theorem for maximum
of the Gaussian sequence on $[0,n]$ as $n\rightarrow\infty.$ Theorem
\ref{Mittal Ylvisaker} easily follows from theirs result, see \cite{book} or
\cite{lectures} for details, definitions, etc. Necessary and sufficient
conditions on $r(k)$ for (\ref{poisson}) are also given in these books.

The main question we considered here is how to approximate the point process
$\mathfrak{B}_{u,n}(B)$ by a family of Poisson processes $\mathfrak{P}%
_{np}(B)$ as level $u$ tends to infinity slower than in (\ref{norm0},
\ref{norm}), that is if $np(u)\rightarrow\infty$ as $u,n\rightarrow\infty$ for
any $B$ of positive measure (length). To this end, using one of the Prokhorov
celebrated theorems, \cite{prokhorov}, on the distance in variation between
Bernoulli and Poisson distributions, we consider the limit behavior of
distributions of random variables $\mathfrak{B}_{u,n}(B)$ with such behaviours
of $u,$ that is, for various medium tending $u$ to infty. First we consider
independent $X(k)$s, and then, using comparison techniques for Gaussian
distributions, pass to dependent ones.

\section{A sequence of independent Gaussian variables}

Assume first that $r(k)=0$ for all $k>0.$ Repeat in the above notations and
conditions a result of Yu. V. Prokhorov on approximations of binomial
distributions. In \cite{prokhorov} approximation qualities of binomial
distribution by both Poisson and normal distributions is considered. Here we
are interesting only in Poisson approximation therefore we formulate only
corresponding result from \cite{prokhorov}

Denote by $|B|,$ the measure (length) of $B,$ and
\begin{equation}
\ \rho_{u,n}(B)=\sum_{k=0}^{n|B|}|P(\mathfrak{B}_{u,n}(B)=k)-P(\mathfrak{P}%
_{np}(B)=k)|, \label{Pr1}%
\end{equation}
$B\in\mathcal{B}$, the distance in variation between the distributions.

\begin{theorem}
\label{Th2} (Theorem 2, \cite{prokhorov}). Assume that $r(k)=0$ for all $k>0.$
Then for all $u\geq0$ and $B\in\mathcal{B}$, $|B|>0,$%
\begin{equation}
\rho_{u,n}(B)=\lambda_{1}p(u)+p(u)O\left(  \min(1,(np(u)^{-1/2}\right)
,\ n,u\rightarrow\infty,\ \label{th1}%
\end{equation}
with $\lambda_{1}=\sqrt{2}/\sqrt{\pi e}=0.483...\ .$
\end{theorem}

From this theorem immediately follows estimation of distance in variance
between point processes $\mathfrak{B}_{u,n}(\cdot)$ and $\mathfrak{P}%
_{np}(\cdot).$

\begin{corollary}
\label{invar0} If $r(k)=0$ for all $k>0,$then
\begin{equation}
\sup_{B\in\mathcal{B}}|\mathfrak{B}_{u,n}(B)-\mathfrak{P}_{np}(B)|\leq Cp(u).
\label{invar-0}%
\end{equation}

\end{corollary}

Indeed, $O(\cdot)$ in (\ref{th1}) is equal to $O(1)$ and does not depend on
$|B|.$

We are interesting here in approximation of $\mathfrak{B}_{u,n}(B)$ for levels
$u$ less than standard one (\ref{norm}). Namely, starting with the above
Prokhorov Theorems we consider levels $u$ tending to infinity with $n$, but
$p(u)n\rightarrow\infty$ as $n\rightarrow\infty.$

Consider an example.

\begin{example}
\label{power scale}\textbf{Power scale. }Let $p(u)n^{a}\rightarrow
c\in(0,\infty).$ For $a=1$ we have Poisson Limit Theorem even for depending
$X(k)$s, that is Theorem \ref{Mittal Ylvisaker}. For $a=0,$ we have Bernoulli
Theorem. A generalization (normal approximation) for depended $X(k)$s one can
find in \cite{book}, \cite{lectures}. Taking in (\ref{norm}) $n^{a},$ $a>0,$
instead of $n,$ we get that%
\begin{equation}
u=u(a)=\sqrt{2a\log n}-\frac{\frac{1}{2}\log\log n+\log(c\sqrt{a\pi/2})}%
{\sqrt{2a\log n}}+O(1/\log n),\ n\rightarrow\infty. \label{example}%
\end{equation}

\end{example}

In \cite{prokhorov} normal approximation (Bernoulli Theorem) for
$\mathfrak{B}_{u,n}(B)$ is considered as well and both the approximations are
compared. In future publications, using the other results from that paper, we
consider normal approximation combined with Poisson one for the number of
level exceedances as well.

\subsection{Thinning. Clusters.}

Consider the following scheme of Poisson approximation of the point process
$\mathfrak{B}_{u,n}(\cdot)$ when $np(u)\rightarrow\infty,$ $n,u\rightarrow
\infty.$ Let an integer $l$ be depended of $n,$ $l=l(n)$ in such a way that
for some $\lambda\in(0,\infty),$
\begin{equation}
\lim_{u,n\rightarrow\infty}(1-p(u))^{l(n)}p(u)=\lambda. \label{tau1}%
\end{equation}
Introduce a thinned point process of $l$\emph{-points of exceedances},
\begin{equation}
\mathfrak{B}_{u,n,l}(B):=\sum_{k\in nB}\mathbf{I}\{\max_{i=1,...,l}X(k-i)\leq
u,X(k)>u\},\ B\in\mathcal{B}. \label{thinned}%
\end{equation}
We call it \emph{cluster center process}, see \cite{Daley V-J}, Section 6.3.

Taking logarithm in (\ref{tau1}) and using Thaylor we get that
\begin{equation}
l(n)=\frac{\log np(u)-\log(\lambda+o(1))}{p(u)+o(1)}=\frac{\log np(u)}%
{p(u)}(1+o(1)),\ \ n,u\rightarrow\infty. \label{tau3}%
\end{equation}
Denote $n_{1}:=l(n)/\log np,$ so that $p(u)n_{1}\rightarrow1$ as
$n,u\rightarrow\infty.$ Denote also
\[
...y_{k}<y_{k+1}<...,
\]
points of $\mathfrak{B}_{u,n,l}(B)$ and associate a point $y_{k}$ with the
\emph{component processes}, that is, \emph{clusters},%
\[
\mathfrak{B}_{u,n_{1}}(B|y_{k})=\sum_{j\in n_{1}B}\mathbf{I}%
\{X(j)>u\}\mathbf{I}\{[y_{k},y_{k+1})\},\ B\in\mathcal{B},\ k\in\mathbb{Z}.
\]
That is we consider independent groups of points of $\mathfrak{B}_{u,n}(B)$
located between $l(n)$-points, say, $l(n)$\emph{-packs}, in another
normalization. Thus we write for any $B\in\mathcal{B},$%
\begin{equation}
\mathfrak{B}_{u,n}(B)=\int_{\mathbb{R}}\mathfrak{B}_{u,n_{1}}(B|y)\mathfrak{B}%
_{u,n,l}(dy)=\sum_{y_{k}\in\mathfrak{B}_{u,n,l}(B)}\mathfrak{B}_{u,n_{1}%
}(B|y_{k})<\infty. \label{cluster PP}%
\end{equation}

Two notmalizations are used above. The first one, $B\rightarrow nB,$ under its
natural (\textquotedblleft correct\textquotedblright) choice (\ref{norm0})
gives convergence to the Poisson process. But sometimes such the choice of
normalization does not capture the entire time interval of interest, therefore
we need $np(u)\rightarrow\infty.$ Therefore first we divide all points of
interest into clusters by $l$-points, and consider points within the clusters
using a corresponding natural normalization $B\rightarrow n_{1}B$,
\textquotedblleft under a magnifying glass\textquotedblright. Then we shall
see that both the number of clusters is approximately Poisson and the number
of points in the cluster is also approximately Poisson, but in a different normalization.

Another situation when such the clusterisation is natural is a case of strong
dependence, that is, $r(k)$ is close to one for several first $k=1,2,...$.
This is especially evident for Gaussian stationary processes in continuous
time $X(t)$ with non-smooth trajectories and the corresponding sequences
$X(k\Delta)$ with $\Delta\rightarrow0$ as $u\rightarrow\infty$, when the
points $t$ of exceeding of high levels appear in the form of short groups
(packs), the number of which is asymptotically Poisson. See details in
\cite{lectures}, Lecture 17 and \cite{book}, Section 4. Let us note that
sometimes normal approximation of the clusters can be more natural.

\subsection{ Convergence %in variation
 to cluster Poisson process.}

Now introduce cluster Poisson process. Cluster center process is a Poisson
point process $\mathfrak{P}_{\lambda}(B),$ $B\in\mathcal{B}$, with intensity
$\lambda,$ see (\ref{tau1}), so that by Theorem \ref{Th2} for any
$B\in\mathcal{B}$, taking $(1-p(u))^{l(n)}p(u)$ instead of $p(u),$ for some
constant $C,$ similarly to Corollary \ref{invar0},
\begin{align}
\sup_{B\in\mathcal{B}}|\mathfrak{B}_{u,n,l}(B)-\mathfrak{P}_{np(u,l)}(B)|  &
\leq Cp(u,l)\ \label{invar-1}\\
\text{with }p(u,l)  &  :=(1-p(u))^{l(n)}p(u).\nonumber
\end{align}
As well, denoting
\[
...z_{k}<z_{k+1}<...,
\]
points of $\mathfrak{P}_{np(u,l)}$, associate a point $z_{k}$ with the Poisson
clusters $\mathfrak{P}_{n_{1}p(u)}(B|z_{k}),$ $k\in\mathbb{Z}$, independent
Poisson point processes with equal intensities $n_{1}p(u).$ Write, for any
$B\in\mathcal{B},$%
\begin{equation}
\mathfrak{P}_{np}(B)=\int_{\mathbb{R}}\mathfrak{P}_{n_{1}p(u)}%
(B|z)\mathfrak{P}_{np(u,l)}(dz)=\sum_{z_{k}\in\mathfrak{P}_{np(u,l)}%
(B)}\mathfrak{P}_{n_{1}p(u)}(B|z_{k})<\infty, \label{cluster PPP}%
\end{equation}
the corresponding cluster Poisson process.

Since all $X(k)$ are independent, the proof of the following theorem obviously
follows from Poisson Limit Theorem and Kallenberg theorem, \cite{kallenberg}.

\begin{theorem}
\label{weak conv} Let $r(k)=0$ for all $k>0.$ The cluster point process
$\mathfrak{B}_{u,n}(B),$ $B\in\mathcal{B}$, (\ref{cluster PP}) converges
weakly as $n,u\rightarrow\infty$ with $np(u)\rightarrow\infty$ to the cluster
Poisson process $\mathfrak{P}_{np}(B)$, $B\in\mathcal{B}$, (\ref{cluster PPP}).
\end{theorem}

Moreover, from Corollary (\ref{invar0}) it follows

\begin{theorem}
\label{var conv} Let $r(k)=0$ for all $k>0.$ The cluster point process
$\mathfrak{B}_{u,n}(B),$ $B\in\mathcal{B}$, (\ref{cluster PP}) converges in
varition as $n,u\rightarrow\infty$ with $np(u)\rightarrow\infty$ to the
cluster Poisson process $\mathfrak{P}_{np}(B)$, $B\in\mathcal{B}$,
(\ref{cluster PPP}).
\end{theorem}

Let us continue Example \ref{power scale}. From (\ref{tau3}) it follows that
\begin{align*}
l(n)  &  =\frac{(a-1)\log(cn/(\lambda+o(1)))}{\log(1-p(u_{n}^{a}))}\\
&  =\frac{1-a}{c}n^{a}\log\frac{cn}{\lambda}\left(  1+o\left(  \frac{1}{\log
n}\right)  \right)  ,\ n\rightarrow\infty.
\end{align*}
Remark that $l(n)=0$ for $a=1.$ Remark as well that $\lambda$ increases up to
infinity as $l$ decreases up to zero, what seems quite naturally.

\section{Dependent Gaussian variables.}

Now turn to Gaussian stationary sequence $X(k)$. We intend to use the
comparison technique for Gaussian distributions, see \cite{berman},
\cite{lectures}, \cite{book}. The following statement is the Corollary 2.3.1,
\cite{lectures}, which is a generalization of Berman inequality,
\cite{berman}, \cite{lectures}, \cite{book}. For any sequence of real numbers
$x_{k},$ $k\in\mathbb{Z}$, denote by $\mathcal{A}_{u},$ algebra of sets
generated by sets $\{x_{k}>u\},$ $k\in\mathbb{Z}$. Denote for shortness by
$\mathbf{X}=\{X(k),k=1,...,n\mathbb{\}}$ and by $\mathbf{X}_{0}=\{X_{0}%
(k),k=1,...,n\mathbb{\}},$where $X_{0}(k)$ are Gaussian independent standard
variables. Let (\ref{logCond}) be fulfilled but in contrast with Theorem
\ref{Mittal Ylvisaker} assume that $np\rightarrow\infty$ as $u\rightarrow
\infty.$ Consider now how fast may $np$ tend to infinity so that Theorem
\ref{var conv} assertion still holds. Similarly to proof of Theorem
\ref{Mittal Ylvisaker}, the main tool is the following comparison inequality.

\begin{proposition}
\label{comparison} For any $A\in\mathcal{A}_{u}$ and any $u,$%
\begin{equation}
|P(\mathbf{X}\in A)-P(\mathbf{X}_{0}\in A)|\leq\frac{1}{\pi}\sum_{k=1}%
^{n}\frac{(n-k)|r(k)|}{\sqrt{1-r^{2}(k)}}\exp\left(  -\frac{u^{2}}%
{1+r(k)}\right)  . \label{3.*}%
\end{equation}

\end{proposition}

Now we derive bounds for $u$ to have that the right hand part of this
inequality still tend to zero as $n\rightarrow\infty.$ Notice that since we
interesting in case $np\rightarrow\infty,$ from (\ref{norm}) and proof of
Mittal result, \cite{mittal ylv}, see also Theorem 3.7, \cite{lectures}, it
can be seen that
\[
\limsup_{n\rightarrow\infty}\frac{u}{\sqrt{2\log n}}\leq1.
\]
Denote $\rho(k):=\sup_{l\geq k}|r(l)|.$

\begin{lemma}
\label{compLemma}Assume that
\begin{equation}
r(k)k^{1-\rho(1)}\rightarrow0\ \text{as }k\rightarrow\infty; \label{powerCond}%
\end{equation}
and
\begin{equation}
\liminf_{n\rightarrow\infty}\frac{u}{\sqrt{2\log n}}>\sqrt{1-\rho(1)}.
\label{lowerbound}%
\end{equation}
Then the sum in (\ref{3.*}) tends to zero as $n\rightarrow\infty.$
\end{lemma}

\textbf{Proof. }(TO CHECK!) Notice that from (\ref{powerCond}) it follows that
$\rho(1)<1.$ Denote $\gamma=1-\rho(1).$ Take some $\alpha\in(0,\gamma)$ and
break the sum in (\ref{3.*}) in two parts, till $[n^{\alpha}]$ and from
$[n^{\alpha}]+1$ till $n.$ We have for the first part,
\begin{align*}
&  \sum_{k=1}^{[n^{\alpha}]}(n-k)|r(k)|\int_{0}^{1}\frac{1}{\sqrt{1-h^{2}%
r^{2}(k)}}\exp\left(  -\frac{u^{2}}{1+hr(k)}\right)  dh\\
&  \leq nn^{\alpha}\rho(1)\frac{1}{\sqrt{1-\rho^{2}(1)}}\exp\left(
-\frac{u^{2}}{1+\rho(1)}\right)  .
\end{align*}
By taking logarithm, we get, that the right part tends to zero as
$n\rightarrow\infty$ if
\[
u^{2}-(1+\rho(1)(1+\alpha)\log n\ \rightarrow\infty
\]
as $n\rightarrow\infty.$ Since $\alpha$ can be chosen arbitrarily, we can say
that there exist $\alpha\in(0,\gamma)$ such that the first part of the sum
(\ref{3.*}) tends to zero as $n\rightarrow\infty$ if
\[
\liminf_{n\rightarrow\infty}\frac{u}{\sqrt{2\log n}}>\sqrt{\frac{1+\rho(1)}%
{2}}.
\]

For the second part we have,
\begin{align*}
&  \sum_{k=[n^{\alpha}]+1}^{n}(n-k)|r(k)|\int_{0}^{1}\frac{1}{\sqrt
{1-h^{2}r^{2}(k)}}\exp\left(  -\frac{u^{2}}{1+hr(k)}\right)  dh\\
&  \leq\frac{n}{\sqrt{1-\rho^{2}(1)}}\sum_{k=[n^{\alpha}]+1}^{n}%
|r(k)|\exp\left(  -\frac{u^{2}}{1+|r(k)|}\right) \\
&  \leq\frac{n}{\sqrt{1-\rho^{2}(1)}}\exp\left(  -\frac{u^{2}}{1+\rho
([n^{\alpha}])}\right)  \sum_{k=[n^{\alpha}]+1}^{n}|r(k)|.
\end{align*}
By condition (\ref{powerCond}), the latter sum for any $\varepsilon>0$ and all
sufficiently large $n$ is at most $\varepsilon n^{2-\gamma}.$ Using this, we
get that the second part of the sum is for the same $\varepsilon$ and $n$ at
most
\[
\frac{\varepsilon n^{2-\gamma}}{\sqrt{1-\rho^{2}(1)}}\exp\left(  -\frac{u^{2}%
}{1+\rho([n^{\alpha}])}\right)  .
\]
Since $\varepsilon$ is arbitrarily small, by taking logarithm, we get that the
second part tends to zero if and only if
\[
\frac{u^{2}}{1+n^{-\alpha\gamma}}-(2-\gamma)\log n\rightarrow\infty
\ \ \text{as }n\rightarrow\infty.
\]
In turn this is followed from the inequality
\[
\liminf_{n\rightarrow\infty}\frac{u}{\sqrt{(2-\gamma)\log n}}>1.
\]
Now just remark that
\[
2-\gamma=1+\rho(1).
\]
Thus the Lemma.

Thus we have proved the following.

\begin{theorem}
\label{final}Let for covariance function of Gaussian sequence $X(k)$ relation
(\ref{powerCond}) be fulfilled. Let $n$ and $u$ tends both to infinity such
that $nP(X(1)>u)\rightarrow\infty$ but (\ref{lowerbound}) be fulfilled. Then
assertion of Theorem \ref{weak conv} is fulfilled in the same notations.
\end{theorem}

\section{Further considerations and extensions.}

\textbf{1. Strong mixing condition.} Remark that condition (\ref{lowerbound})
means that $\sup_{k\geq1}|r(k)|<1/2.$ If $\sup_{k\geq1}|r(k)|\geq1/2$ Lemma
\ref{compLemma} does not work. But one can apply Theorem 3.5, \cite{lectures},
see also Theorem 2.1, \cite{book}, from which it follows that if
\[
\sum_{k=1}^{\infty}k|r(k)|<\infty,
\]
then for any $u$ and any normalization, point process $\eta_{u}(\cdot)$
satisfies Rosenblatt strong mixing condition. Hence assertion of Theorem
\ref{final} can be proved by this approach for any normalizations of the mark
processes $\zeta_{u}^{k}(\cdot).$ It is a subject of future considerations.

\textbf{2. Brown motion clusters. }If some applications require the normalization of
$\zeta_{u}^{k}(\cdot)$ with $pn_{1}\rightarrow\infty,$ one can prove
convergence to a cluster Poisson process with independent Wiener processes
with trends as clusters, or marks.

\textbf{3. Random energy model by Derrida,} \cite{derrida},

For independent $X(k),$ the classical Derrida model of randon energy is
\[
S_{N}(\beta):=\sum_{k=1}^{[2^{N}]}e^{\beta\sqrt{N}X(k)},\ \beta>0.
\]
Standard problems here are study a limit behavior of $S_{N}(\beta)/N$ as
$N\rightarrow\infty,$ in dependence of $\beta$ and limit behavior of
distribution of normed $S_{N}(\beta)$ (limit theorem), as well. Since large
values of the sequence $X(k)$ plays main role, the above results can be apply.
Our approach allows also considering dependent $X(i)$s. Moreover, Prokhorov
theorems may allow to get quality of these limit approximations and even to
derive asymptotic expansions.

\end{document}